\documentclass[12pt]{article}
\usepackage{amsmath,amssymb,amsthm}
\usepackage[mathscr]{eucal}
\usepackage{rotate,graphics,epsfig}
\usepackage{color}
\usepackage{hyperref}

\newtheorem{The}{Theorem}
\newtheorem{Exa}[The]{Example}

\newtheorem{Pro}[The]{Proposition}
\newtheorem{Lem}[The]{Lemma}
\theoremstyle{definition}

\newtheorem{Rem}[The]{Remark}
\numberwithin{equation}{section}
\numberwithin{The}{section}

\newcommand{\be}{\begin{eqnarray}}
\newcommand{\ee}{\end{eqnarray}}

\newcommand{\by}{\begin{eqnarray*}}
\newcommand{\ey}{\end{eqnarray*}}

\newcommand{\bn}{\begin{enumerate}}
\newcommand{\en}{\end{enumerate}}

\newcommand{\bi}{\begin{itemize}}
\newcommand{\ei}{\end{itemize}}

\def\frac#1#2{{#1 \over #2}}

\begin{document}
\title{Tail Dependence of Multivariate Archimedean Copulas}
\author{
Haijun Li\footnote{{\small\texttt{lih@math.wsu.edu}}, Department of Mathematics and Statistics, Washington State University, 
Pullman, WA 99164, U.S.A.}
}
\date{December 2024}
\maketitle

\begin{abstract}

Archimedean copulas generated by Laplace transforms have been extensively studied in the literature, with much of the focus on tail dependence limited only to cases where the Laplace transforms exhibit regular variation with positive tail indices. In this paper, we extend the investigation to include Archimedean copulas associated with both slowly varying and rapidly varying Laplace transforms. We show that tail dependence functions with various tail orders effectively capture the extremal dependence across the entire class of Archimedean copulas, reflecting the full spectrum of tail behaviors exhibited by the underlying Laplace transforms.


\end{abstract}

\section{Introduction}
\label{S1}

Archimedean copulas constitute a fundamental class of copula functions that are widely used in the copula theory and applications (see, e.g., \cite{Nelsen2006, Joe2014, JL11}), and can be best illustrated by using scale mixtures of independently and identically distributed exponential random variables. 
Let $E_1, \dots, E_d$ be independent and exponentially distributed random variables with unit mean, and $R$ be a strictly positive random variable that is independent of $E_i$, $1\le i\le d$. Define:
\begin{equation}
	\label{mix exponential}
	{  X}=(X_1, \dots, X_d) := (RE_1, \dots, RE_d).
\end{equation}
Assume that the Laplace transform of $R^{-1}$ exists and is given by $\varphi(t)$, which is continuous and decreasing with $\varphi(0)=1$ and $\lim_{t \to \infty}\varphi(t)=0$. The marginal survival function of $X_i$ is then described as  $\overline{F}_i(t) = \mathbb{P}(X_i>t)= \mathbb{E}(e^{-t/R}) = \varphi(t)$, $t\ge 0$, which implies that its inverse  $\overline{F}_i^{\,-1}(u_i) = \varphi^{-1}(u_i)$, for $0\le u_i\le 1$, $1\le i\le d$. The joint distribution function of the transformed scale mixture $(\varphi(X_1), \dots, \varphi(X_d))$, after applying $\varphi$ to $X$, is now derived as follows.
\begin{eqnarray}
	{C}(u_1 \dots, u_d) &:=& \mathbb{P}\big(RE_1>\overline{F}_1^{\,-1}(u_1), \dots, RE_d>\overline{F}_d^{\,-1}(u_d)\big)= \mathbb{E}\big(e^{-[\sum_{i=1}^d\varphi^{-1}(u_i)]/R}\big)\nonumber\\
	&=& \varphi\Big(\sum_{i=1}^d\varphi^{-1}(u_i)\Big),~\forall~ (u_1, \dots, u_d)\in [0,1]^d.\label{Archimedean copula}
\end{eqnarray}
The function \eqref{Archimedean copula} is known as an Archimedean copula with {\em generator} $g:=\varphi^{-1}$. In general, a copula $C: [0,1]^d\to [0,1]$ is defined as a distribution function of a random vector $(U_1, \dots, U_d)$ with univariate, uniformly distributed margins $U_i$ on $[0,1]$, $1\le i\le d$; see \cite{Joe97, Joe2014} for detailed discussions on copulas and applications. Correspondingly, the distribution $\widehat{C}$ of the dual vector $(1-U_1, \dots, 1-U_d)$ is called the survival copula. As demonstrated from scale mixture \eqref{mix exponential}, any random vector can be decomposed into univariate margins and a copula, that encodes all the marginally scale-free dependence information, so that marginal behaviors and dependence structure can be studied separately.

Observe that the Laplace transform $\varphi$ of a positive random variable is {completely monotone} on $[0,\infty)$; that is, $\varphi$ has the derivatives of all orders and
\begin{equation}
	\label{complete monotone}
	(-1)^k\frac{d^k\varphi(t)}{dt^k}\ge 0, ~\forall~t\in [0,\infty), ~ k=0, 1, 2, \dots.
\end{equation} 
In fact, since the Laplace transforms of positive random variables coincide with complete monotone functions on $[0, \infty)$ with a value of 1 at 0, the condition \eqref{complete monotone} is also sufficient for \eqref{Archimedean copula} to be a copula, according to Kimberling's theorem (see \cite{Kimberling74}). 
\begin{Pro}\rm
	\label{Kimberling} 
	Let $\varphi: \mathbb{R}_+\to [0,1]$ be continuous and strictly decreasing with boundary conditions $\varphi(0)=1$ and $\lim_{t\to \infty}\varphi(t)=0$. The function \eqref{Archimedean copula} defined on $[0,1]^d$ is a copula for any $d\ge 2$ if and only if $\varphi$ is completely monotone. 
\end{Pro}

It is observed from \eqref{Archimedean copula} that the decay rate of probabilities of the scale mixture $X$ in upper-orthants, or equivalently, the lower-orthant decay rate of Archimedean copula $C$, is determined by right-tail behaviors of Laplace transform $\varphi$. The most commonly examined univariate right-tail pattern is the power-law behavior, also known as regular variation, which is defined as that $\varphi$ is said to be {\em regularly varying at $\infty$}  with tail parameter $-\alpha$, denoted by $\varphi\in \mbox{RV}_{-\alpha}$, if for all $x>0$, 
\begin{equation}
	\label{univ RV}
	\frac{\varphi(tx)}{\varphi(t)}\to x^{-\alpha}, \ \mbox{as}\ t\to \infty.
\end{equation}
It is straightforward that $\varphi\in \mbox{RV}_{-\alpha}$ if and only if $\varphi(x) = \ell(x)x^{-\alpha}$ for some $\ell(x)\in \mbox{RV}_{0}$, which is known as {\em slowly varying at $\infty$} \cite{Resnick07}. 
On the other hand, a function $\varphi(x)$ is said to be {\em regularly varying at $0$} with tail parameter $\beta$ if $\varphi(x^{-1})$ is regularly varying at $\infty$ with tail parameter $-\beta$.

The goal of this paper is to complete a tail dependence theory of  Archimedean copulas with various Laplace transforms $\varphi$, including but not limited to regular variation. In fact, it is well-known that if $\varphi\in \mbox{RV}_{-\alpha}$, where $\alpha>0$, then the lower-orthant decay rate of $C$ can be explicitly derived, see, e.g., \cite{CS09, JLN10}. In Section 2, we complement this result by deriving the lower-orthant decay rate of $C$ for a {slowly varying} Laplace transform $\varphi\in \mbox{RV}_0$. Moreover, in Section 3, we derive the lower-orthant decay rate of an Archimedean copula $C$ when the Laplace transform $\varphi$ is {\em rapidly varying} (in the sense of Omey \cite{Omey2013}). In short, the tail dependence of an Archimedean copula emerges from all three distinct right-tail decay patterns of its underlying Laplace transform $\varphi$: slow variation, regular variation, and rapid variation.

The lower-orthant decay rate of a copula $C$ can be evaluated by using the $k$-th order tail dependence function (\cite{JLN10, HJ11}), denoted by $b(w; k)$, $k\ge 1$, that is defined as the non-zero limiting function 
\begin{equation}
\label{k tail limit}
b(w; k)= \lim_{u\downarrow 0}\frac{C(uw_1, \dots, uw_d)}{u^k\,\ell(u)}=\lim_{u\downarrow 0}\frac{\mathbb{P}\big(\cap_{i=1}^d\{U_i\le uw_i\}\big)}{u^k\,\ell(u)}, 
\end{equation} 
for $w=(w_1, \dots, w_d)\in (0,\infty)^d$, 
whenever the limit exists, for some function $\ell(\cdot)$, that is slowing varying at $0$; i.e., $\ell(t^{-1})\in \mbox{RV}_0$. Observe that the tail order $k$ satisfies $1\le k\le d$ for a non-trivial $k$-th order tail dependence function to exist. We discuss in this paper the lower tail dependence of copula $C$ only, as the upper tail dependence of $C$ can be studied by analyzing the lower tail dependence of the survival copula $\widehat{C}$.

If $k=1$, the tail dependence function $b(w;1)$, first introduced in \cite{JLN10}, can be trivially extended to $\mathbb{R}^d_+$. The existences of the first-order tail dependence for $C$ and its margins are equivalent to the existence of the first-order exponent function, denoted by $a(w;1)$, that is defined as the non-zero limiting function 
\begin{equation}
	\label{exponent function}
	a(w; 1)= \lim_{u\downarrow 0}\frac{\mathbb{P}\big(\cup_{i=1}^d\{U_i\le uw_i\}\big)}{u\,\ell(u)}, \ w=(w_1, \dots, w_d)\in \mathbb{R}^d_+, 
\end{equation} 
whenever the limit exists, for some function $\ell(\cdot)$, that is slowing varying at $0$. The exponent function of a copula, together with univariate, regularly varying margins, are shown in \cite{Li2009, LS2009, LH13} to be equivalent to multivariate regular varying distributions. Multivariate regular variation describes multivariate power-law behaviors, and has been extensively studied in the literature; see, e.g., \cite{Resnick07, MS01}.

If $k>1$, the tail dependence function $b(w;k)$ describes the higher order tail dependence of multivariate distributions that is often emerged from hidden regular variation within a subcone of $\mathbb{R}^d\backslash \{0\}$ or multivariate rapid variation (\cite{HJ11, HJL12, LH14}). For example, if $X$ is distributed according to a multivariate normal distribution with mean zero and $d\times d$ equicorrelation matrix, having pair-wise correlation coefficient $0\le \rho<1$, then its first order tail dependence parameter $b((1,\dots, 1); 1) = 0$, but the higher tail order $k = d/[1+(d-1)\rho]$ (see \cite{HJ11}), satisfying that $1<k\le d$, that reveals the higher order tail dependence emerged from the multivariate normal distribution with exponentially decayed joint tails (also known as multivariate rapid variation). In this paper, we demonstrate that various tail orders also exist for Archimedean copulas with slowly varying, regularly varying, and rapidly varying Laplace transforms $\phi$.

Throughout this paper, the right-tail equivalence $f(x)\sim g(x)$ means that two non-zero functions satisfy that $f(x)/g(x)\to 1$ as $x\to a$, $a\in [0,\infty]$.

\section{Tail dependence of Archimedean copulas with regularly varying Laplace transforms}
\label{S2}

Since an Archimedean copula is generated from the Laplace transform $\varphi$ (or generator $\varphi^{-1}$) of a positive random variable, its joint tail decay structure depends on the tail behaviors of $\varphi$, as was shown in the case of regular variation in \cite{genest1989, CS09, JLN10}. We provide the proofs here for compeleteness.

\begin{The}\rm
	\label{RV} Let $C({u}; \varphi)=\phi(\sum_{i=1}^d\varphi^{-1}(u_i))$
	be an Archimedean copula, $u=(u_1, \dots, u_d) \in [0,1]^d$, with Laplace transform $\varphi$.  
	\begin{enumerate}
		\item If the Laplace transform $\varphi$ is
		regularly varying at $\infty$ with tail parameter $-\alpha<0$, then the
		lower tail dependence function of $C$ is given by 
		\begin{equation}
			\label{lower tail-RV}
			b(w;1) =\Bigl(\sum_{j=1}^d w_j^{-1/\alpha}\Bigr)^{-\alpha}, \ w = (w_1, \dots, w_d)\in \mathbb{R}_+^d. 
		\end{equation}
	\item If the
	inverse Laplace transform (or generator) $\varphi^{-1}$ is regularly varying at
	1, or $\varphi^{-1}(1-1/x)$ is regularly
	varying at $\infty$ with tail parameter $-\beta<0$, then the upper exponent function of $C$ is given by
		\begin{equation}
		\label{upper tail-RV}
		a(w;1) =\Bigl(\sum_{j=1}^d w_j^{\beta}\Bigr)^{1/\beta}, \ w = (w_1, \dots, w_d)\in \mathbb{R}_+^d. 
	\end{equation}
	\end{enumerate}
\end{The}

\noindent
{\sl Proof.} (1)
Using the stochastic representation \eqref{mix exponential}, the tail dependence function \eqref{lower tail-RV} is a natural consequence from the regular variation of $\varphi$.  In fact, it follows from Proposition 2.6 of \cite{Resnick07} that  $\overline{F}_i^{\,-1}(u_i) = \varphi^{-1}(u_i)=\ell(u_i)u_i^{-1/\alpha}$ is regularly varying at $0$ with tail parameter $-1/\alpha$, where $\ell(\cdot)$ is slowing verying at $0$. Therefore, for $w = (w_1, \dots, w_d)>0$, 
\[\varphi^{-1}(uw_i)=\ell(uw_i)(uw_i)^{-1/\alpha} \sim \ell(u)u^{-1/\alpha}w_i^{-1/\alpha}=\varphi^{-1}(u)w_i^{-1/\alpha},\ 1\le i\le d. 
\]
Since $\varphi(\cdot)$ is continuous, we have by \eqref{Archimedean copula}, 
\[b(w;1) = \lim_{u\downarrow 0}\frac{C(uw_1, \dots, uw_d)}{u} = \lim_{u\downarrow 0}\frac{\varphi\Big(\varphi^{-1}(u)\sum_{i=1}^dw_i^{-1/\alpha}\Big)}{\varphi(\varphi^{-1}(u))}= \Bigl(\sum_{j=1}^d w_j^{-1/\alpha}\Bigr)^{-\alpha},
\]
which trivially holds if some of $w_i$s are zero, and hence \eqref{lower tail-RV} holds 
for all $w = (w_1, \dots, w_d)\in \mathbb{R}_+^d$. 

(2) 
Similarly, the upper tail dependence of Archimedean copula $C({ u};
\varphi)=\phi\bigl(\sum_{i=1}^d\varphi^{-1}(u_i)\bigr)$ emerges when the
inverse Laplace transform (or generator) $\varphi^{-1}$ is regularly varying at
1; that is, $\varphi^{-1}(1-1/x)$ is regularly
varying at $\infty$ with tail parameter $-\beta<0$ (see \cite{genest1989} for the bivariate case). For all $w = (w_1, \dots, w_d)\in \mathbb{R}_+^d$, 
\[\varphi^{-1}(1-uw_i)=\ell(uw_i)(uw_i)^{\beta} \sim \ell(u)u^{\beta}w_i^{\beta}=\varphi^{-1}(1-u)w_i^{\beta},\ 1\le i\le d, 
\]
implying that
\[\varphi(\varphi^{-1}(1-u)w_i)\sim 1-uw_i^{1/\beta}, \ 1\le i\le d. 
\]
Observe that 
\[\mathbb{P}(U_i>1-uw_i,\ \exists\ i)=1-\mathbb{P}(U_i\le 1-uw_i,\ \forall\ i)=1-\varphi\big(\sum_{i=1}^d \varphi^{-1}(1-uw_i)\big)
\]
\[\sim 1-\varphi\big(\varphi^{-1}(1-u)\sum_{i=1}^dw_i^{\beta}\big)\sim 1-\Big(1-u\big(\sum_{i=1}^dw_i^\beta\big)^{1/\beta}\Big)\sim u\big(\sum_{i=1}^dw_i^\beta\big)^{1/\beta}. 
\]
Since $\varphi(\cdot)$ is continuous, we have, 
\[a(w;1) = \lim_{u\downarrow 0}\frac{\mathbb{P}(U_i>1-uw_i,\ \exists\ i)}{u} = \Bigl(\sum_{j=1}^d w_j^{\beta}\Bigr)^{1/\beta},
\]
and \eqref{upper tail-RV} holds. 
\hfill $\Box$

\begin{Rem}
	\label{Archimedean-r-1}
\begin{enumerate}
	\item If the Laplace transform $\varphi$ is
	regularly varying at $\infty$ with tail parameter $-\alpha<0$, then $X$ in \eqref{mix exponential} is multivariate regularly varying at $\infty$, with upper intensity measure $\mu_\infty(\cdot)$; that is, for any fixed norm  $||\cdot||$ on $\mathbb{R}^d_+$,
	\begin{equation}
	\label{MRV at infinity}
	\frac{\mathbb{P}(X\in tB)}{\mathbb{P}(||X||>t)}\to \mu_\infty(B),\ \mbox{as}\ t\to \infty, 
	\end{equation} 
for all the relative compact subsets $B\subset \overline{\mathbb{R}}^d_+\backslash\{0\}$ bounded away from $0$, that is $\mu_\infty$-continuous in the sense  that $\mu_\infty(\partial B)=0$. The intensity measure $\mu_\infty(\cdot)$ is known to enjoy the group invariance $\mu_\infty(tB) = t^{-\alpha}\mu_\infty(B)$, $t>0$, leading to the semi-parametric representation for multivariate extremes that is useful in statistical analysis \cite{Resnick07}. 
	Since an Archimedean copula $C$ can be viewed as the survival copula of $X$ in \eqref{mix exponential}, the intensity measure $\mu_\infty(\cdot)$ yields the lower tail dependence function of $C$, as presented in Theorem \ref{RV} (1). 
		\item It is straightforward to see that the
	inverse Laplace transform $\varphi^{-1}$ is regularly varying at
	1 if and only if $F_i(u) = 1 -\varphi(u)$ is regularly varying at $0$ with tail parameter $1/\beta$. If $\varphi^{-1}$ is regularly varying at
	1, then $X$ in \eqref{mix exponential} is multivariate regularly varying at $0$, with lower intensity measure $\mu_0(\cdot)$; that is,  
		\begin{equation}
		\label{MRV at zero}
		\frac{\mathbb{P}(X\in uB)}{\mathbb{P}(||X||\le u)}\to \mu_0(B)
	\end{equation} 
	for all the relative compact subsets $B\subset \overline{\mathbb{R}}^d_+\backslash\{+\infty\}$, that is $\mu_0$-continuous in the sense  that $\mu_0(\partial B)=0$ \cite{MS01}. Similar to Remark \ref{Archimedean-r-1} (1), \eqref{MRV at zero}
	yields the upper tail dependence function of $C$, as presented in Theorem \ref{RV} (2). 
\end{enumerate}
\end{Rem}

We focus on multivariate maximums in this paper, and hence it follows from Remark \ref{Archimedean-r-1} that we discuss only the tail dependence function $b(\cdot;k)$ for Archimedean copulas with Laplace transform $\varphi$. 
The proof of Theorem \ref{RV}, however, excludes the case that $\alpha=0$, when the Laplace transform $\varphi$ is slowly varying at infinity. To obtain the tail dependence in this case, the following lemma is needed. 

\begin{Lem}\rm (Elez and Djur\v{c}i\'c \cite{Elez2013})
	\label{Inv-SV} Let $g: [0,\infty)\to (0,\infty)$ be slowing varying with inverse $g^{-1}$ and $\lim_{x\to \infty}g(x) = \infty$. Then the inverse $g^{-1}$ is rapidly varying in the sense of de Haan \cite{Haan1970}; that is,
	\[\lim_{y\to \infty}\frac{g^{-1}(\lambda y)}{g^{-1}(y)}= \infty,\ \ \forall\ \lambda>1.
	\]
\end{Lem}

\begin{The}\rm
	\label{SV} Let $C({u}; \varphi)=\phi(\sum_{i=1}^d\varphi^{-1}(u_i))$
	be an Archimedean copula, $u=(u_1, \dots, u_d) \in [0,1]^d$,  where the Laplace transform $\varphi$ is
	slowly varying at $\infty$. The
	lower tail dependence function of $C$ is given by 
	\begin{equation}
		\label{lower tail-SV}
		b(w;1) =\min\{w_1, \dots, w_d\}, \ w = (w_1, \dots, w_d)\in \mathbb{R}_+^d. 
	\end{equation}
\end{The}

\noindent
{\sl Proof.} Since $\varphi(x)$ is slowly varying and decreasing to zero,  $1/\varphi(x)$ is slowly varying and increasing to $\infty$. Observe that 
\[\left(1/\varphi\right)^{-1}(y) = \varphi^{-1}(1/y),\ y>0.
\]
Therefore, by Lemma \ref{Inv-SV}, $\varphi^{-1}(1/y)$ is rapidly varying in the sense of de Haan; that is, for any $\lambda>1$, 
\[\frac{\varphi^{-1}\left(1/(\lambda y)\right)}{\varphi^{-1}(1/y)}\to \infty,\ \ y\to \infty. 
\]
Let $u=1/y$, and thus, for any $\lambda>1$, 
\[\frac{\varphi^{-1}\left(u/\lambda\right)}{\varphi^{-1}(u)}\to \infty,\ \ u\to 0,
\]
which is equivalent to 
\begin{equation}
		\label{inv-SV0}
\frac{\varphi^{-1}\left(\lambda u\right)}{\varphi^{-1}(u)}\to 0,\ \ u\to 0, \forall\ \lambda>1. 
\end{equation}

For any fixed $w = (w_1, \dots, w_d)>0$, let $m=\min\{w_1, \dots, w_d\}>0$ and $\lambda_i = w_i/m\ge 1$, $1\le i\le d$. Consider 
\[\frac{\varphi\left(\sum_{i=1}^d\varphi^{-1}(uw_i)\right)}{u}= m\frac{\varphi\left(\sum_{i=1}^d\varphi^{-1}(\lambda_ium)\right)}{um} = m\frac{\varphi\left(\sum_{i=1}^d\frac{\varphi^{-1}(\lambda_ium)}{\varphi^{-1}(um)}\varphi^{-1}(um)\right)}{um}. 
\]
It follows from \eqref{inv-SV0} that as $u\to 0$, 
\[\sum_{i=1}^d\frac{\varphi^{-1}(\lambda_ium)}{\varphi^{-1}(um)}\to c_w
\]
where $c_w=|\{i|w_i=m, 1\le i\le d\}|\ge 1$ is a constant. Hence, for a fixed $0<\epsilon<1$, there exists a small $\delta>0$, such that whenever $0<u<\delta$, 
\[c_w-\epsilon\le \sum_{i=1}^d\frac{\varphi^{-1}(\lambda_ium)}{\varphi^{-1}(um)}\le c_w+\epsilon. 
\]

Since the Laplace transform $\varphi(\cdot)$ is decreasing, 
\[m\frac{\varphi\left((c_w-\epsilon)\varphi^{-1}(um)\right)}{um}\ge \frac{\varphi\left(\sum_{i=1}^d\varphi^{-1}(uw_i)\right)}{u}\ge m\frac{\varphi\left((c_w+\epsilon)\varphi^{-1}(um)\right)}{um}.
\]
When $u\to 0$, $\varphi^{-1}(um)\to \infty$. It then follows from the slow variation of $\varphi$ that as $u\to 0$, 
\[\frac{\varphi\left((c_w+\epsilon)\varphi^{-1}(um)\right)}{um}=\frac{\varphi\left((c_w+\epsilon)\varphi^{-1}(um)\right)}{\varphi(\varphi^{-1}(um))}\to 1
\]
\[\frac{\varphi\left((c_w-\epsilon)\varphi^{-1}(um)\right)}{um}=\frac{\varphi\left((c_w-\epsilon)\varphi^{-1}(um)\right)}{\varphi(\varphi^{-1}(um))}\to 1, 
\]
implying that $b(w;1)=\lim_{u\downarrow 0}\frac{\varphi\left(\sum_{i=1}^d\varphi^{-1}(uw_i)\right)}{u}=m$, for any $w>0$. The tail dependence \eqref{lower tail-SV} trivially holds if some of $w_i$s are zero, and hence \eqref{lower tail-SV} holds for any $w\in \mathbb{R}_+^d$. 
\hfill $\Box$

\begin{Rem}
\begin{enumerate}
				\item In stochastic representation \eqref{mix exponential}, if $R$ is regularly varying at $\infty$ with tail parameter  $-\alpha\le 0$, then, by Breiman's Theorem (see \cite{Breiman1965}), the margin $X_i$ is regularly varying with tail parameter  $-\alpha\le 0$. In particular, if $R$ is slowly varying at $\infty$, then $\varphi\in \mbox{RV}_{0}$. 
	\item Observe from \eqref{mix exponential} that $\varphi(x)$ is the marginal survival function of $X=(X_1, \dots, X_d)$. The fact that $\varphi\in \mbox{RV}_{-\alpha}$, $\alpha\ge 0$, indicates that $X$ in \eqref{mix exponential} is multivariate regularly varying in the sense of \eqref{MRV at infinity}, including multivariate slowly varying, with a simple mixture structure, from which tail dependence emerges. 
	\end{enumerate}
\end{Rem}

\begin{Exa}\rm 
	It is well-known that the Laplace transforms of positive random variables coincide with complete monotone functions on $[0, \infty)$. When $d=2$, in particular, the class of Laplace transforms consists of non-increasing, non-negative convex functions on $[0, \infty)$. 
	
	Let $\varphi(x) = 1/\log(x+e)$, $x\ge 0$. The function is non-negative, decreasing, convex and slowly varying as $x\to \infty$, and hence $\varphi(x)$ is the Laplace transform of a positive random variable, that is slowly varying. The inverse $\varphi^{-1}(u) = e^{1/u}-e$, $0\le u\le 1$, and the straightforward calculation leads to 
	\[b(w_1, w_2;1) = \lim_{u\downarrow 0}\frac{\varphi(\varphi^{-1}(uw_1)+\varphi^{-1}(uw_2))}{u}= \min\{w_1, w_2\}, 
	\]
	for $(w_1,w_2)\in \mathbb{R}^2_+$. 
\end{Exa}

\section{Tail dependence of Archimedean copulas with rapidly varying Laplace transforms}
\label{S3}

The higher order tail dependence naturally emerges from an Archimedean copula with rapidly varying Laplace transform $\varphi$. Let $\Gamma_\alpha(g)$ denote the class of measurable functions for which there exists a measurable and positive function $g$, such that
\begin{equation}
	\label{gamma class}
	\lim_{x\to \infty}\frac{f(t+xg(t))}{f(t)}=e^{-\alpha x},\ \ \forall x\in \mathbb{R},\ \alpha\ge 0.
\end{equation} 
The $\Gamma_\alpha(g)$-class is studied extensively in \cite{Omey2013}, where $\alpha$ can be any real number. We assume that $\alpha\ge 0$ because we focus on non-increasing functions in $\Gamma_\alpha(g)$ only in this paper. 
It is shown in \cite{Omey2013} that the limit \eqref{gamma class} holds locally uniformly, where $g$ is self-neglecting in the sense that $g(t)/t\to 0$ and 
\begin{equation}
	\label{self-neglecting}
	\lim_{x\to \infty}\frac{g(t+xg(t))}{g(t)}=1,\ \ \forall x\in \mathbb{R},
\end{equation} 
holds locally uniformly. The local uniform convergences yield the representations for the functions in  $\Gamma_\alpha(g)$, $\alpha\ge 0$. In particular, the following two representations from \cite{Omey2013} are needed.
\begin{enumerate}
	\item If $g$ is self-neglecting, then $g(x) = D(x)W(x)$, where $D(x)\to c>0$, and
\begin{equation}
	\label{self-neglecting rep}
 W(x)=\exp\Big\{\int_{x_0}^x\frac{\epsilon^*(z)}{g(z)} dz\Big\}
\end{equation} 
such that its derivative converges to zero. That is, $g$ can be taken as a function whose derivative converges to zero. Note that $\epsilon^*(z)$ can be negative, whereas $g(z)>0$, $z\in \mathbb{R}$. 
\item If $g$ is self-neglecting and $f$ has a non-increasing derivative $f'$, then $g$ and $-f/f'$ are tail equivalent. That is, if $f\in \Gamma_\alpha(g)$ is the survival function of a random variable, then the self-neglecting function $g$ can be taken as the reciprocal of the hazard rate.
\end{enumerate}

\begin{Lem}\rm 
	\label{ultimately IFR} If $\varphi\in \Gamma_\alpha(g)$, where $g$ is ultimately monotone (i.e., monotone for all sufficiently large $t$), then the self-neglecting function $g$ satisfies that the limit $g(t)/g(\lambda t)$, as $t\to \infty$, exists (real numbers or positive infinity) for any $\lambda>1$. 
\end{Lem}

\noindent
{\sl Proof.}
We prove the ultimately decreasing case only, and the ultimately increasing case is similar. 
 It follows from \eqref{self-neglecting rep} that $\epsilon^*(z)$ is ultimately negative (i.e., negative for all sufficiently large $z$). Observe that 
\[\frac{g(t)}{g(\lambda t)}\sim \exp\Big\{\int_{\lambda t}^t\frac{\epsilon^*(z)}{g(z)} dz\Big\},\ \mbox{as}\ t\to \infty, 
\]
is ultimately increasing, and therefore, the limit $g(t)/g(\lambda t)$ exists for any $\lambda>1$. 
\hfill $\Box$

\begin{Rem}
	The Laplace transform $\varphi(t)\in \Gamma_\alpha(g)$ is the survival function of a positive random variable $X_i$ (see \eqref{mix exponential}), and $g$ can be taken as the reciprocal of the failure rate of $X_i$. Therefore, $g(t)$ is ultimately monotone can be interpreted as the failure rate of $X_i$ being ultimately monotone.
	\end{Rem}

Since $\varphi(x)$ is the marginal survival function of $X=(X_1, \dots, X_d)$, the rapid variation of the Laplace transform $\varphi\in \Gamma_{\alpha}(g)$ indicates that $X$ in \eqref{mix exponential} is multivariate rapidly varying in some sense, with a simple scale mixture structure, which yields tail dependence, as the following result shows.  

\begin{The}\rm
	\label{Rapid V} Let $C({u}; \varphi)=\phi(\sum_{i=1}^d\varphi^{-1}(u_i))$
	be an Archimedean copula, $u=(u_1, \dots, u_d) \in [0,1]^d$,  where the Laplace transform $\varphi\in \Gamma_{\alpha}(g)$ is 
	rapidly varying at $\infty$, $\alpha>0$, where $g$ is ultimately decreasing. If, in addition, the Laplace transform $\varphi$ satisfies that for some $1\le k\le d$, 
	\begin{equation}
		\label{ratio rapid}
		\lim_{t \to \infty}\frac{\varphi(td)}{\varphi^k(t)}= \tau >0,
	\end{equation} 
	then the
	lower tail dependence function of $C$ is given by 
	\begin{equation}
		\label{tail-Rapid}
		b(w;k) =\tau\prod_{i=1}^dw_i^{k/d}, \ w = (w_1, \dots, w_d)\in \mathbb{R}_+^d. 
	\end{equation}
\end{The}

\noindent
{\sl Proof.}
Since $\varphi(t)$ is the marginal survival function of $X_i$ in the stochastic representation \eqref{mix exponential} and $g$ is the reciprocal of its hazard rate, one can write
\[\varphi(x) = \exp\Big\{-\int_0^x\frac{1}{g(z)}dz\Big\}, x\ge 0. 
\]
It follows from Lemma \ref{ultimately IFR} that the limit $g(t)/g( td)$ exists, and thus, applying L'H\^opital's rule on \eqref{ratio rapid} yields
\[\tau = \lim_{t \to \infty}\frac{\varphi(td)d/g(td)}{\varphi^k(t)k/g(t)}=\tau \lim_{t \to \infty}\frac{d}{k}\frac{g(t)}{g(td)}
\]
resulting in that $\lim_{t \to \infty}g(t)/g(td) = k/d$.

Let $u=\varphi(t)$, and thus $u$ is small if and only if $t$ is large. Furthermore let $e^{-\alpha x_i}=w_i$, $1\le i\le d$. Observe that
\[\varphi\big(t+g(t)(\log w_i^{-1/\alpha})\big)=\varphi(t+x_ig(t))\sim \varphi(t)e^{-\alpha x_i}=uw_i, \ (x_1, \dots, x_d)\in \mathbb{R}^d_+,
\]
imply that
\[td+g(t)\sum_{i=1}^d\log w_i^{-1/\alpha}\sim \sum_{i=1}^d\varphi^{-1}(uw_i).
\]
Take $\varphi$ on both sides, and therefore, 
\[\varphi\Big(\sum_{i=1}^d\varphi^{-1}(uw_i)\Big)\sim \varphi\Big(td+g(t)\sum_{i=1}^d\log w_i^{-1/\alpha}\Big)\sim \varphi\Big(td+kg(td)/d\sum_{i=1}^d\log w_i^{-1/\alpha}\Big),
\]
which, in turns, is tail equivalent to 
$\varphi(td)\exp\Big\{-\alpha\sum_{i=1}^d\log w_i^{-k/d\alpha}\Big\}=\varphi(\varphi^{-1}(u)d)\prod_{i=1}^dw_i^{k/d}
$. Hence, it follows from \eqref{ratio rapid} that 
\[
b(w;k) = \lim_{u\downarrow 0}\frac{\varphi\Big(\sum_{i=1}^d\varphi^{-1}(uw_i)\Big)}{u^k}=\tau\lim_{u\downarrow 0}\frac{\varphi\Big(\sum_{i=1}^d\varphi^{-1}(uw_i)\Big)}{\varphi(\varphi^{-1}(u)d)}=\tau\prod_{i=1}^dw_i^{k/d},
\]
for $1\le k\le d$, as desired. 
\hfill $\Box$

\begin{Rem}
	\label{Archimedean-r-5}
	\begin{enumerate}
		\item Huang obtained this result in \cite{Huang2020} under an additional condition that the self-neglecting function $g$ satisfies that $g(t)/g(td)$ converges to a constant, as $t\to \infty$. Such a restriction is not necessary in our new proof, due to Lemma \ref{ultimately IFR}. 
		\item If $g$ is ultimately decreasing, then the failure rate of $X_i$, $1\le i\le d$, is ultimately increasing, and hence $X_i$, $1\le i\le d$, has the light right tail. On the other hand, if $g$ is ultimately increasing, then the failure rate of $X_i$, $1\le i\le d$, is ultimately decreasing, and tail behaviors are more complex with possible heavier tails.  
		\item Most distributions are defined in terms of densities, and the tail density approach is used in \cite{LH14, JL19, Li2021} to study  the multivariate rapid variation and higher order tail dependence functions. As in the case of multivariate regular variation, if a density is multivariate rapidly varying, then the corresponding distribution is multivariate rapidly varying. In contrast, however, the densities of Archimedean copulas often do not enjoy compact explicit forms. 
	\end{enumerate}	
\end{Rem}

It should be mentioned that the tail order $k$ depends on $\alpha$ through \eqref{ratio rapid} implicitly, in contrast to the regular variation case of Theorem \ref{RV}, where the tail dependence function depends on $\alpha$ explicitly. The assumption \eqref{ratio rapid}, however, is mild and satisfied in various applications. 
Several examples are given below to illustrate Theorem \ref{Rapid V}, and all these Archimedean copulas and the corresponding Laplace transforms can be found in \cite{Joe97}.

\begin{Exa}\rm 
Consider a bivariate Frank copula $C(u_1,u_2)$ with Laplace transformation of $\varphi(t) = -\frac{\log (1-(1-e^{-\theta})e^{-t})}{\theta}$, $\theta>0$. Since $\varphi(t)\sim \frac{1-e^{-\theta}}{\theta}e^{-t}$, then, let $g(t) \equiv 1$, 
		\begin{eqnarray}
	\varphi(t+g(t)x) &\sim & \frac{1-e^{-\theta}}{\theta} e^{-(t+x)}      
	\sim \varphi(t)e^{- x}. \nonumber
\end{eqnarray}
That is, $\varphi\in \Gamma_1(g(t))$, and $g(t)/g(2t)=1$ and $\alpha =1$. Observe that
\[
\frac{\varphi(2t)}{\varphi(t)^2} \sim  \frac{(\frac{1-e^{-\theta}}{\theta})(e^{-t})^2}{(\frac{1-e^{-\theta}}{\theta} e^{-t})^2} = \frac{\theta}{1-e^{-\theta}}=\tau.  
\]
According to Theorem \ref{Rapid V}, its lower tail dependence $b(w_1, w_2; 2) = \frac{\theta}{1-e^{-\theta}}w_1w_2$, $w_1, w_2\ge 0$.  
\end{Exa}

\begin{Exa}\rm 
	Consider a bivariate B5 copula $C(u_1,u_2)$ with Laplace transformation of $\varphi(t) = 1 - (1-e^{-t})^{1/\theta}, \theta\ge 1$. Since  $\varphi(t) \sim (1/\theta) e^{-t},$ then, let $g(t)\equiv 1$, 
	\[\varphi(t+1x) \sim  (1/\theta) e^{-t} e^{-x}  
	=\varphi(t)e^{-\alpha x},
	\]
	where $\alpha = 1$. 
	That is, $\varphi\in \Gamma_1(g(t))$ and $g(t)/g(2t)=1$. Observe that 
	\[\frac{\varphi(td)}{\varphi(t)^2} \sim \frac{\frac{1}{\theta}(e^{-t})^{2}}{(\frac{1}{\theta}e^{-t})^{2}} = \theta = \tau.
	\]
	According to Theorem \ref{Rapid V}, its lower tail dependence $b(w_1, w_2; 2) = \theta w_1w_2$, $w_1, w_2\ge 0$.  
\end{Exa}

\begin{Exa}\rm 
	Consider a two-dimensional Archimedean copula $C(u_1,u_2)$ with Laplace transformation of $\varphi(t) = [(1-\theta)e^{-t}/(1-\theta e^{-t})]^\alpha, 0\le \theta< 1$ and $\alpha>0$. This is the Laplace transform of a negative binomial
	distribution. Observe that
	\[\varphi(t)\sim (1-\theta)^\alpha e^{-\alpha t},\ 0\le \theta< 1, \alpha>0. 
	\]
	\[\varphi(t+1x)\sim (1-\theta)^\alpha e^{-\alpha (t+x)}\sim \varphi(t) e^{-\alpha x}
	\]
	where $g(t)\equiv 1$ is self-neglecting, leading to $\varphi\in \Gamma_{\alpha}(g(t))$, and $g(t)/g(2t)=1$. 
		Since,
	\begin{eqnarray}
		\frac{\varphi(td)}{\varphi(t)^2} &\sim & \frac{(1-\theta)^\alpha(e^{-\alpha t2})}{((1-\theta)^\alpha e^{-\alpha t})^{2}} = (1-\theta)^{-\alpha} = \tau,
	\end{eqnarray}
the lower tail dependence $b(w_1, w_2; 2) = (1-\theta)^{-\alpha}w_1w_2$, $w_1, w_2\ge 0$. 
\end{Exa}

\begin{Exa}\rm 
	Consider a two-dimensional  Gumbel copula, $C(u_1,u_2)$ with Laplace transformation  $\varphi(t) = e^{-t^{1/\theta}},\ \theta\ge 1$. 
	Let $g(t)=t^{1-1/\theta}$, where $\theta\ge 1$, which is non-decreasing. It is straightforward to verify that $g(t)$ is self-neglecting, and 
	\[\frac{g(t)}{g(2t)}= 2^{-1+1/\theta}.
	\]
	Observe that as $t\to \infty$, 
	\[\frac{\varphi(t+y\,g(t))}{\varphi(t)}\sim \frac{e^{-t^{1/\theta}\big[1+y\,t^{-1/\theta}/\theta\big]}}{e^{-t^{1/\theta}}}\to e^{-y/\theta}
	\]
	and thus $\varphi(t)\in \Gamma_{1/\theta}(g(t))$. Furthermore, we also have
	\[\frac{\varphi(2t)}{\varphi(t)^{2^{1/\theta}}}= e^{-t^{1/\theta}2^{1/\theta}+t^{1/\theta}2^{1/\theta}}=1,
	\] and thus the lower tail dependence function is given by 
	\[b(w_1,w_2; 2^{1/\theta}) = \big(w_1w_2\big)^{2^{-1+1/\theta}},
	\]
where 	$w_1 \ge 0$,  and  $w_2 \ge 0$.
\end{Exa}

\section{Concluding remarks}
\label{S4}

The extremal dependence of a random vector 
$X=(X_1, \dots, X_d)$ is influenced by, and thus can vary with, the tail behaviors of its marginal distributions. In contrast, tail dependence, as introduced in \cite{JLN10}, depends solely on the copula of $X$, independent of its marginal distributions. This inherent extremal dependence captures a universal characteristic present in all multivariate distributions, including those with multivariate slowly varying, regularly varying, and rapidly varying joint tails. In this paper, we illustrate this universality via the tail dependence function with tail order $k$ (\cite{JLN10, HJ11}) for random samples, where the marginal distributions, as characterized by Laplace transform $\varphi$,  exhibit diverse tail decay patterns ranging from slow to regular to rapid variation.

Due to their simple scale-mixing structure, Archimedean copulas are typically specified by their distribution functions, whereas most multivariate distributions are defined through their densities \cite{Resnick07, BE07}. Tail densities of order 
$k$, introduced in \cite{LW2013, LH14}, have been used to analyze the universal characteristics of extremal dependence via the copula density approach. Additionally, \cite{JL19, Li2021} derived the tail dependence of skew-elliptical distributions, capturing the universal features of extremal dependence in both regularly varying and rapidly varying skew-elliptical distributions. More generally, the density-based approach can be used for analyzing tail dependence through distributional densities.


\end{document}